\documentclass{amsart}

\usepackage{amssymb,color}
\usepackage{amsfonts}
\usepackage{amsmath}
\usepackage{euscript}
\usepackage{enumerate}
\usepackage{graphics}

\usepackage{hyperref}
\usepackage{pdfsync}
\newtheorem*{Theorem1'}{Theorem 1'}

\theoremstyle{definition}

\theoremstyle{remark}

\newcommand \R{{\mathbb R}}
\newcommand \Z{{\mathbb Z}}
\newcommand \Q{{\mathbb Q}}
\newcommand \C{{\mathbb C}}
\newcommand \N{{\mathbb N}}
\renewcommand \S{\mathcal{S}}

\renewcommand \P{\mathcal{P}}
\def\a{{\alpha}}
\def\b{{\beta}}
\renewcommand\chi{{c}}

\begin{document}

\title[On the degree of repeated radical extensions]{On the degree of repeated radical extensions}

\author{Fernando Szechtman}
\address{Department of Mathematics and Statistics, University of Regina, Canada}
\email{fernando.szechtman@gmail.com}
\thanks{This research was partially supported by an NSERC grant}



\dedicatory{This paper is dedicated to Natalio H. Guersenzvaig.}


\begin{abstract} We answer a question posed by Mordell in 1953, in the case repeated radical extensions, and find necessary and sufficient conditions for $[F[\sqrt[m_1]{N_1},\dots,\sqrt[m_\ell]{N_\ell}]:F]=m_1\cdots m_\ell$,
where $F$ is an arbitrary field of characteristic not dividing any $m_i$.
\end{abstract}

\maketitle










\section{Introduction}

We fix throughout a unique factorization domain $D$ with field of fractions $F$, allowing for
the possibility that $D=F$, and write $\chi(F)$ for the characteristic of $F$. We also fix $\ell\in\N$, $m_1,\dots,m_\ell\in\N$, $m=\mathrm{lcm}\{m_1,\dots,m_\ell\}$, and $N_1,\dots,N_\ell\in D$. A prime means a prime positive integer.

In this paper, we give a necessary and sufficient condition for
\begin{equation}
\label{extension}
[F[\sqrt[m_1]{N_1},\dots,\sqrt[m_\ell]{N_\ell}]:F]=m_1\cdots m_\ell,
\end{equation}
assuming only $\chi(F)\nmid m$. This settles a problem posed by Mordell \cite{M} in 1953, in the case of repeated radical extensions.

The degrees of repeated radical extensions have been studied by several authors, including Hasse~\cite{H},
Besicovitch \cite{B}, Mordell \cite{M}, Siegel \cite{S}, Richards \cite{Ri}, Ursell \cite{U}, Jian-Ping~\cite{J}, Albu~\cite{A},
and Carr and O'Sullivan \cite{CS}.

The question of when $[F[\sqrt[n]{a}]:F]=n$ was solved by Vahlen \cite{V} in 1895 if $F=\Q$,
Capelli \cite{C} in 1897 if $F$ has characteristic~0, and R\'edei \cite[Theorem 428]{R} in 1959 in general.

\medskip

\noindent{\bf Irreducibility Criterion (C).} The polynomial $X^n-a\in F[X]$ is irreducible if and only if $a\notin F^p$
for every prime factor $p$ of~$n$, and if $4|n$ then $a\notin -4F^4$.

\medskip

In particular, if $-a\notin F^2$ and $a\notin F^p$ for every prime factor $p$ of~$n$, then $X^n-a$ is irreducible. 
The special case of (C) when $n$ is prime is due to Abel; a very simple proof of this case can be found in \cite[Theorem 427]{R}.

\medskip

Provided $F$ contains a primitive $m$th root of unity, Hasse \cite{H} showed that (\ref{extension}) holds if and only~if
\begin{equation}
\label{Hasse}
(\sqrt[m_1]{N_1})^{a_1}\times\cdots\times (\sqrt[m_\ell]{N_\ell})^{a_\ell}\in F, a_i\geq 0,\text{ only when }m_1|a_1,\dots,m_\ell|a_\ell.
\end{equation}

Later, Besicovitch \cite{B} proved (\ref{extension}) assuming: $D=\Z$; each $N_i$ is positive and has a  prime factor that divides it only once and does not divide any other $N_j$; each $\sqrt[m_i]{N_i}$ is positive and real (the $m_1\cdots m_\ell$ embeddings of
$\Q[\sqrt[m_1]{N_1},\dots, \sqrt[m_\ell]{N_\ell}]$ into $\C$ then yield (\ref{extension}) for all other $m_i$th roots of the $N_i$). The special case of Besicovitch's result when $m_1=\dots=m_\ell$ and every $N_i$ is prime appears in Richards~\cite{Ri} (for the more elementary case
$m_1=\dots=m_\ell=2$, see \cite{F,Ro}). Assuming that $N_1,\dots,N_\ell$ are pairwise relatively prime, Ursell \cite{U} obtained
a variation of Besicovitch's theorem.

Mordell \cite{M} combined and extended the results of Hasse and Besicovitch, and proved (\ref{extension}) assuming (\ref{Hasse}), and that $F$ contains a primitive $m$th root of unity or that $F$ is a subfield of $\R$ with all $\sqrt[m_1]{N_1},\dots, \sqrt[m_\ell]{N_\ell}$ real. In the latter case, all $N_i$ such that $m_i$ is even must be positive. Siegel~\cite{S} gave a theoretical
description of the value of $[F[\sqrt[m_1]{N_1},\dots,\sqrt[m_\ell]{N_\ell}]:F]$ under Mordell's condition that $F$ be a subfield of $\R$ with all $\sqrt[m_i]{N_i}$ real.  Albu \cite{A} extended the work of Mordell and Siegel to the case when $F$ contains a primitive $m$th root of unity or all $m$th roots of unity in $F[\sqrt[m_1]{N_1},\dots,\sqrt[m_\ell]{N_\ell}]$ belong to $\{1,-1\}$. Under these weaker assumptions, (\ref{extension}) is still shown in \cite{A} to be a consequence of (\ref{Hasse}). It is worth noting that, except for (\ref{Hasse}), none of the aforementioned conditions
are necessary for (\ref{extension}) to hold. A different approach was taken by Jian-Ping~\cite{J}, using valuation theory,
when $F$ is an algebraic number field; he succeeded in avoiding any assumptions on roots of unity
and proved a more general version of (\ref{extension}), applicable to repeated extensions via Eisenstein polynomials, not just binomials.
 Nevertheless, Jian-Ping's hypotheses are also unnecessary for (\ref{extension}) to hold. Indeed, when the ring of integers of
$F$ is a UFD, each $N_i$ is forced to have an irreducible factor that divides it only once and does not divide any other $N_j$. More recently, Carr and O'Sullivan \cite{CO} proved a fairly general result on the linear independence of roots and reproved Mordell's theorem as an application.

Set $J=\{1,\dots,\ell\}$, $\P=\{p\,|\, p\text{ is a prime factor of }m\}$, and for each $i\in J$ and $p\in\P$ let $m_i(p)$ be the $p$-part of $m_i$, so that $m_i(p)=p^{n_i}$, where $n_i\geq 0$, $p^{n_i}|m_i$ and $p^{n_i+1}\nmid m_i$.  It is clear that
\begin{equation}
\label{exten}
[F[\sqrt[m_1]{N_1},\dots,\sqrt[m_\ell]{N_\ell}]:F]=m_1\cdots m_\ell\Leftrightarrow [F[\sqrt[m_1(p)]{N_1},\dots,\sqrt[m_\ell(p)]{N_\ell}]:F]=m_1(p)\cdots m_\ell(p)
\end{equation}
for all $p\in\P$. We are thus reduced to study the case when each $m_i=m_i(p)$ for a fixed prime $p$. We split this case in two subcases
depending on the parity of $p$. For each prime $p$, we set
$$
\S_{p}=\{N_1,N_1^{e_1}N_2,N_1^{e_1}N_2^{e_2}N_3,\dots, N_1^{e_1}\cdots N_{\ell-1}^{e_{\ell-1}}N_\ell, 0\leq e_i<p\}.
$$
In particular, $\S_{2}$ consists of all $N_1^{e_1}\cdots N_{\ell}^{e_{\ell}}$ such that $e_i\in\{0,1\}$ and 
$(e_1,\dots,e_\ell)\neq (0,\dots,0)$.

\bigskip

\noindent{\bf Theorem A.}  Let $n_1,\dots,n_\ell\in\N$, $p$ an odd prime such that $\chi(F)\neq p$, and suppose $m_i=p^{n_i}$ for all 
$i\in J$. Then (\ref{extension}) holds if and only if $\S_{p}\cap D^p=\emptyset$.

\bigskip

The well-known example $[\Q[\sqrt[4]{-1},\sqrt[4]{2}]:\Q]=8$ shows that Theorem A fails if $p=2$. The above criteria impose
general conditions that disallow this example. Close examination of numerous pathological cases led us to the exact conditions required
when $p=2$. We say that 
$(N_1,\dots,N_\ell)$ is $2$-defective
if the following two conditions hold: $\S_2\cap D^2=\emptyset$ but $\S_2\cap (-D^2)\neq\emptyset$ 
(this readily implies that $|\S_2\cap (-D^2)|=1$,
as shown in Lemma 5); if $-d^2=M=N_1^{f_1}\cdots N_\ell^{f_\ell}$ is the only element of $\S_2\cap (-D^2)$, where $d\in D$, $0\leq f_i<2$, and
$M^\sharp=\{i\in J\,|\, f_i=1\}$ is nonempty (since $\S_2\cap D^2=\emptyset$, the exponents $f_i$ are uniquely determined by $M$, whence $M^\sharp$ is well-defined), then $4|m_i$ for all $i\in M^\sharp$, and if $i\in M^\sharp$, then
\begin{equation}
\label{spe}
\pm 2d\underset{j\neq i}\Pi N_j^{e_j}\in D^2\text{ for some choice of }0\leq e_j<1.
\end{equation}
Since $-M=N_1^{f_1}\cdots N_\ell^{f_\ell}\in D^2$, the outcome of (\ref{spe}) is independent of the actual choice of $i\in M^\sharp$.

\bigskip

\noindent{\bf Theorem B.}  Let $n_1,\dots,n_\ell\in\N$ and suppose that $\chi(F)\neq 2$ and $m_i=2^{n_i}$ for all $i\in J$. Then (\ref{extension}) holds if and only if $\S_2\cap D^2=\emptyset$ and $(N_1,\dots,N_\ell)$ is not 2-defective.

\bigskip

Combining  (\ref{exten}) with Theorems A and B we immediately obtain a general criterion for (\ref{extension}). This requires additional
notation. For each $p\in\P$, we set $J(p)=\{i\,|\, i\in J\text{ and }p|m_i\}$, and write
$$
J(p)=\{i(p,1),\dots,i(p,\ell(p))\},\quad i(p,1)<\dots< i(p,\ell(p)),
$$
$$
\S(p)=\{N_{i(p,1)},N_{i(p,1)}^{e_1}N_{i(p,2)},N_{i(p,1)}^{e_1}N_{i(p,2)}^{e_2}N_{i(p,3)},\dots,N_{i(p,1)}^{e_1}\cdots
N_{i(p,\ell(p)-1)}^{e_{\ell(p)-1}}N_{i(p,\ell(p))}, 0\leq e_j<p\}.
$$

\bigskip

\noindent{\bf Theorem C.} Suppose that $\chi(F)\neq m$. Then (\ref{extension}) holds if and only if $\S(p)\cap D^p=\emptyset$
for every $p\in\P$ and, if $2\in\P$, then $(N_{i(2,1)},\dots,N_{i(2,\ell(2))})$ is not 2-defective.

\bigskip

The next example illustrates the use of Theorems B and C, and lies outside of the scope of the aforementioned criteria.

\bigskip

\noindent{\bf Example 1.} Suppose that $-1\notin F^2$ and each of $A,B,C\in D$ has an irreducible factor that divides it only
once and does not divide any of the two other elements. Then
$$
[F[\sqrt[m_1]{AB},\sqrt[m_2]{BC},\sqrt[m_3]{-CA}]:F]=m_1m_2m_3
$$
if and only if at least one of $m_1,m_2,m_3$ is not divisible by 4 or none of $\pm 2A,\pm 2B,\pm 2C\in D^2$.

\bigskip

As we are dealing with a classical and basic problem, we purposely resort to elementary and complete arguments
in order maximize the potential readership of our solution.

\section{Lemmata}

Given a nonzero $a\in F$ we write $\langle a\rangle$ for the subgroup of $F^\times$ generated by $a$.

\medskip

\noindent{\bf Lemma 1.} Let $p$ be a prime such that $\chi(F)\neq p$ and suppose $b_1,\dots,b_n\in F$ are nonzero. Then
\begin{equation}\label{le1}
F[\sqrt[p]{b_1},\dots,\sqrt[p]{b_n}]^p\cap F=F^p\langle b_1,\dots,b_n\rangle.
\end{equation}

\medskip

\noindent{\sc Proof.} Let $M,N\in F$ be nonzero. We claim that if $\sqrt[p]{M}\in F[\sqrt[p]{N}]$ then
$M\in F^p\langle N\rangle$. This is clear if $N\in F^p$ so we assume $N\notin F^p$.

Set $K=F[\zeta]$, where $\zeta$ is a primitive $p$th root of unity. Then $K/F$ is a Galois extension
with Galois group isomorphic to a subgroup of $(\Z/p\Z)^\times$. In particular, $[K:F]$ divides $(p-1)$.
It follows that $K^p\cap F=F^p$. Indeed, suppose $a\in F$ and $\a\in K$ satisfies $\a^p=a$. Since $[F[\a]:F]$
divides $[K:F]$, it also divides $p-1$. As $p\nmid (p-1)$, $X^p-a\in F[X]$ is reducible, whence $a\in F^p$ by (C).

By assumption, $N\notin F^p$. Thus $N\notin K^p$ as indicated above, so $X^p-N\in K[X]$ is irreducible by (C).
Thus $\{1,\sqrt[p]{N},\dots,\sqrt[p]{N^{p-1}}\}$ is a $K$-basis of $K[\sqrt[p]{N}]$. By assumption, $\sqrt[p]{M}\in K[\sqrt[p]{N}]$, so
\begin{equation}
\label{ko}
\sqrt[p]{M}=a_0+a_1\sqrt[p]{N}+\cdots+a_{p-1}\sqrt[p]{N^{p-1}},\quad a_i\in K.
\end{equation}
Note that $K[\sqrt[p]{N}]/K$ is a Galois extension with cyclic Galois group $\langle\sigma\rangle$, where
$\sigma(\sqrt[p]{N})=\zeta \sqrt[p]{N}$. Since $\sqrt[p]{M}$ is a root of $X^p-M$, we must have
$\sigma(\sqrt[p]{M})=\zeta^i\sqrt[p]{M}$ for some $0\leq i<p$. Applying $\sigma$ to (\ref{ko}), we obtain
$$
\zeta^i\sqrt[p]{M}=a_0+a_1\zeta\sqrt[p]{N}+\cdots+a_{p-1}\zeta^{p-1}\sqrt[p]{N^{p-1}}.
$$
On the other hand, multiplying (\ref{ko}) by $\zeta^i$ yields
$$
\zeta^i\sqrt[p]{M}=a_0\zeta^i +a_1\zeta^i\sqrt[p]{N}+\cdots+a_{p-1}\zeta^i\sqrt[p]{N^{p-1}}.
$$
From the $K$-linear independence of $1,\sqrt[p]{N},\dots,\sqrt[p]{N^{p-1}}$ we infer
that $a_j=0$ for all $j\neq i$. Thus
$$
\sqrt[p]{M}=a \sqrt[p]{N^i},\quad a\in K,
$$
whence $M=a^p N^i$. Thus $M/N^{-i}\in K^p\cap F=F^p$, so $M\in F^p\langle N\rangle$.

\medskip

By above, $F[\sqrt[p]{b_1}]^p\cap F=F^p\langle b_1\rangle$. Suppose $n>1$ and $F[\sqrt[p]{b_1},\dots,\sqrt[p]{b_{n-1}}]^p\cap F=F^p\langle b_1,\dots,b_{n-1}\rangle$. Then
$$
\begin{aligned}
F[\sqrt[p]{b_1},\dots,\sqrt[p]{b_{n-1}},\sqrt[p]{b_n}]^p\cap F &=
(F[\sqrt[p]{b_1},\dots,\sqrt[p]{b_{n-1}}][\sqrt[p]{b_n}])^p\cap F[\sqrt[p]{b_1},\dots,\sqrt[p]{b_{n-1}}]\cap F\\
&=
F[\sqrt[p]{b_1},\dots,\sqrt[p]{b_{n-1}}]^p\langle b_{n}\rangle\cap F.
\end{aligned}
$$
Let $\a\in F[\sqrt[p]{b_1},\dots,\sqrt[p]{b_{n-1}}]^p\langle b_n\rangle\cap F$. Then $\a\in F$ and
$\a b_{n}^{i}\in F[\sqrt[p]{b_1},\dots,\sqrt[p]{b_{n-1}}]^p\cap F$ for some $i\in\Z$. Thus $\a b_{n}^{i}\in F^p\langle b_1,\dots,b_{n-1}\rangle$
and therefore $\a\in F^p\langle b_1,\dots,b_{n-1},b_{n}\rangle$.\qed

\bigskip

\noindent{\bf Lemma 2.} Suppose $-1\notin F^2$ and $\pm a\notin F^2$. Then for any $n\in\N$, we have $\sqrt{-1}\notin F[\sqrt[2^n]{a}]$.

\bigskip

\noindent{\sc Proof.} We show by induction that $\sqrt{-1}\notin F[\sqrt[2^n]{a}]$ and $\pm \sqrt[2^n]{a}\notin F[\sqrt[2^n]{a}]^2$.
The fact that $\sqrt{-1}\notin F[\sqrt{a}]$ follows from Lemma 1. Suppose, if possible, that $\pm \sqrt{a}=z^2$, where $z\in F[\sqrt{a}]$. Then $z=x+y\sqrt{a}$,
where $x,y\in F$, so that $\pm\sqrt{a}=x^2+ay^2+2xy\sqrt{a}$. It follows that $x^2+ay^2=0$. As $-a\notin F^2$, we infer $x=y=0$,
a contradiction.

Assume we have shown that $\sqrt{-1}\notin F[\sqrt[2^n]{a}]$ and $\pm \sqrt[2^n]{a}\notin F[\sqrt[2^n]{a}]^2$ for some $n\in\N$.
Since $\pm a\notin F^2$, (C) implies $[F[\sqrt[2^n]{a}]:F]=2^n$ and  $[F[\sqrt[2^{n+1}]{a}]:F]=2^{n+1}$. But $\sqrt{-1}\notin F[\sqrt[2^n]{a}]$,
so $F[\sqrt[2^{n+1}]{a}]=F[\sqrt[2^{n}]{a},\sqrt{-1}]$. Thus $\sqrt[2^{n+1}]{a}=\a+\b\sqrt{-1}$ for unique $\a,\b\in F[\sqrt[2^{n}]{a}]$. Squaring, we get $\sqrt[2^{n}]{a}=\a^2-\b^2+2\a\b\sqrt{-1}$, which implies $\a\b=0$ and $\a^2-\b^2=\sqrt[2^n]{a}$, a contradiction.\qed

\bigskip

\noindent{\bf Lemma 3.} Suppose $\chi(F)\neq 2$, let $n\in\N$, and set $K=F[\zeta]$, where $\zeta$ is a primitive $2^n$th root of unity. If $n\leq 2$ or $-1\in F^2$, then $G=\mathrm{Gal}(K/F)$ is cyclic.

\bigskip

\noindent{\sc Proof.} We have an embedding $\Psi:G\to (\Z/2^n\Z)^\times$, $\sigma\to [s]$, where
$\sigma(\zeta)=\zeta^s$. This settles the case $n\leq 2$. Assume henceforth that $n\geq 3$.
We have $(\Z/2^n\Z)^\times=\langle a,b\rangle$, where $a=[5]$, $b=[-1]$ and $\langle a\rangle\cap\langle b\rangle$
is trivial \cite[Chapter VI]{Vi}.  By hypothesis, $-1=\alpha^2$, where $\alpha\in F\cap\langle\zeta\rangle$. Suppose, if possible,
that $b\in\Psi(G)$, say $b=\Psi(\sigma)$. Then $\sigma(\alpha)=\alpha^{-1}=-\alpha$, since $\alpha$ is a power of $\zeta$,
and $\sigma(\a)=\a$, since $\a\in F$. This contradiction shows that $b\notin \Psi(G)$. Now any subgroup $S$ of $\langle a,b\rangle$
that does not contain $b$ must be cyclic (if $S$ is not trivial, it is generated by $a^i$ or $a^i b$, where $i$ is the smallest
positive integer such that an element of this type is in $S$). Thus $G$ is cyclic.\qed

\bigskip

\noindent{\bf Lemma 4.} Let $n\in\N$, $p$ an odd prime such that $\chi(F)\neq p$, and set $K=F[\zeta]$, where $\zeta$ is a primitive $p^n$th root of unity. Then $\mathrm{Gal}(K/F)$ is cyclic.

\bigskip

\noindent{\sc Proof.} $\mathrm{Gal}(K/F)$ is isomorphic to a subgroup of $(\Z/p^n\Z)^\times$,
which is a cyclic group.\qed

\bigskip

\noindent{\bf Lemma 5.}  Suppose $\S_2\cap D^2=\emptyset$. Then $|S_2\cap (-D^2)|\leq 1$, with $|S_2\cap (-D^2)|=0$ if $-1\in F^2$.

\bigskip

\noindent{\sc Proof.} Suppose $M\neq N$ are in $\S_2\cap (-D^2)$. Then $MN\in D^2$ and $MN=e^2P$, where $P\in\S_2$ and $e\in D$. Thus $P\in D^2$, against $\S_2\cap D^2=\emptyset$. If $-1\in F^2$ then $-D^2=D^2$, so $\S_2\cap (-D^2)=\emptyset$.\qed

\bigskip

\noindent{\bf Lemma 6.} Suppose $\chi(F)\neq 2$, let $n\in\N$ and set $K=F[\zeta]$, where $\zeta$ is a primitive $2^n$th root of unity. Assume $\S_2\cap D^2=\emptyset$. Then
$|\S_2\cap K^2|\in \{0,1,3\}$. Moreover, if $|\S_2\cap K^2|=3$ then one of the elements of
$\S_2\cap K^2$ is in $\S_2\cap (-D^2)$, and we have $-1\notin F^2$, $n\geq 3$.

\bigskip

\noindent{\sc Proof.} Suppose $M\neq N$ are in $\S_2\cap K^2$. Then $MN\in K^2$ and $MN=e^2P$, where $P\in\S_2$ and $e\in D$, so $P\in D^2$.
Lemma 1 implies that $F[\sqrt{M}],F[\sqrt{N}],F[\sqrt{P}]$ are distinct intermediate
subfields of $K/F$ of degree 2. In particular, $\mathrm{Gal}(K/F)$ is not cyclic. Now  $\mathrm{Gal}(K/F)$ is isomorphic to a subgroup of $(\Z/2^n\Z)^\times$, so $n\geq 3$ and $(\Z/2^n\Z)^\times\cong (\Z/2^{n-2}\Z)\times (\Z/2\Z)$. Any subgroup of
$(\Z/2^{n-2}\Z)\times (\Z/2\Z)$ has at most 3 subgroups of index 2, so the Galois correspondence implies that any intermediate
subfield of $K/F$ of degree 2 must be equal to one of $F[\sqrt{M}],F[\sqrt{N}],F[\sqrt{P}]$. Lemma 1 readily implies that no
element from $\S_2$ different from $M,N,P$ is in $K^2$. By Lemma 3, $-1\notin F^2$, so
$F[\sqrt{-1}]$ must be equal to one of $F[\sqrt{M}],F[\sqrt{N}],F[\sqrt{P}]$, and Lemma 1 implies that one of $M,N,P$ is in $-D^2$.\qed

\bigskip

\noindent{\bf Lemma 7.}  Let $n_1,\dots,n_\ell\in\N$ and $p$ a prime such that $\chi(F)\neq p$ and $m_i=p^{n_i}$ for all $i\in J$.
Let $K=F[\zeta]$, where $\zeta$ is a primitive $m$th root of unity, $m=\mathrm{lcm}\{m_1,\dots,m_\ell\}$.
Suppose that $\S_p\cap K^p=\emptyset$. Then $[K[\sqrt[m_1]{N_1},\dots,\sqrt[m_\ell]{N_\ell}]:K]=m_1\cdots m_\ell$.

\bigskip

\noindent{\sc Proof.} By assumption $N_1\notin K^p$. Moreover, if $4|m_1$ then $-1\in K^2$ and therefore $-N_1\notin K^2$.
It follows from (C) that $[K[\sqrt[m_1]{N_1}]:K]=m_1$.
Suppose $[K[\sqrt[m_1]{N_1},\dots,\sqrt[m_i]{N_i}]:K]=m_1\cdots m_i$ for some $1\leq i<\ell$.

Assume, if possible,
that $\sqrt[p]{N_{i+1}}\in K[\sqrt[m_1]{N_1},\dots,\sqrt[m_i]{N_i}]$. Then $K[\sqrt[p]{N_{i+1}}]$
is an intermediate subfield of degree $p$ in the Galois extension $K[\sqrt[m_1]{N_1},\dots,\sqrt[m_i]{N_i}]/K$, with Galois group
$G=\langle \sigma_1,\dots,\sigma_i\rangle$, where
$$
\sigma_k(\sqrt[m_k]{N_k})=\zeta^{m/m_k}\, \sqrt[m_k]{N_k},\; \sigma_k(\sqrt[m_j]{N_j})=\sqrt[m_j]{N_j},\; j\neq k.
$$
Any subgroup of $G$ of index $p$ contains $G^p$, so by the Galois correspondence $K[\sqrt[p]{N_{i+1}}]$ is contained in the fixed field of $G^p$, namely $K[\sqrt[p]{N_1},\dots,\sqrt[p]{N_i}]$. Lemma 1 implies that $N_1^{e_1}\cdots N_i^{e_{i}}N_{i+1}\in K^p$ for some $0\leq e_i<p$,
against $\S_p\cap K^p=\emptyset$. Thus $\sqrt[p]{N_{i+1}}\notin K[\sqrt[m_1]{N_1},\dots,\sqrt[m_i]{N_i}]$.

Assume if possible, that $4|m_{i+1}$ and $\sqrt{-N_{i+1}}\in K[\sqrt[m_1]{N_1},\dots,\sqrt[m_i]{N_i}]$. Then the above argument yields
$N_1^{e_1}\cdots N_i^{e_{i+1}}N_{i+1}\in -K^2$ for some $0\leq i<p$. But $-K^2=K^2$, so $\S_p\cap K^p=\emptyset$ is violated.
This shows $\sqrt{-N_{i+1}}\notin K[\sqrt[m_1]{N_1},\dots,\sqrt[m_i]{N_i}]$ when $4|m_{i+1}$.

We deduce from (C) that $[K[\sqrt[m_1]{N_1},\dots,\sqrt[m_{i+1}]{N_{i+1}}]:K]=m_1\cdots m_{i+1}$.\qed

\bigskip

\noindent{\bf Lemma 8.}  Let $n_1,\dots,n_\ell\in\N$ and $p$ a prime such that $\chi(F)\neq p$ and $m_i=p^{n_i}$ for all $i\in J$.
Let $K=F[\zeta]$, where $\zeta$ is a primitive $m$th root of unity, $m=\mathrm{lcm}\{m_1,\dots,m_\ell\}$. 
Suppose that $\S_p\cap D^p=\emptyset$ and $|\S_p\cap K^p|=1$,
say $M=N_1^{f_1}\cdots N_\ell^{f_\ell}$, where $0\leq f_i<p$, and $M^\sharp=\{i\in J\,|\, f_i=1\}$ is nonempty. For $i\in M^\sharp$,
set $V_i=\{\sqrt[m_1]{N_1},\dots,\sqrt[m_\ell]{N_\ell}\}\setminus \{\sqrt[m_i]{N_i}\}$ and let $m[i]$ be the product of all $m_j$ with 
$j\neq i$.
Then

(a) $[K[V_i]:K]=m[i]$ and $\sqrt[p]{N_i}\notin F[V_i]$ for all $i\in M^\sharp$.

(b) If $p$ is odd or $m_i=2$ for at least one $i\in M^\sharp$, then (\ref{extension}) holds.

(c) If $4|m$ and $\S_2\cap (-D^2)=\emptyset$, then  $\sqrt{-N_i}\notin F[V_i]$ for all $i\in M^\sharp$, so (\ref{extension}) holds.

(d) If $4|m_i$ for all $i\in M^\sharp$ and $M\in \S_2\cap (-D^2)$, say $M=-d^2$ with $d\in D$, then (\ref{extension}) holds if and only if given any $i\in M^\sharp$,
(\ref{spe}) fails.

\bigskip

\noindent{\sc Proof.} Let $i\in M^\sharp$. By Lemma 7, we have $[K[V_i]:K]=m[i]$ and hence $[F[V_i]:F]=m[i]$. Suppose, if possible, that
 $\sqrt[p]{N_i}\in F[V_i]$ and set $Y_i=\{\zeta\}\cup V_i$. Then
  $F[\sqrt[p]{N_i}]$
is an intermediate subfield of degree $p$ in the Galois extension $F[Y_i]/F$,
with Galois group $G=H\rtimes U$, where $H=\langle \sigma_j\,|\, j\neq i\rangle$ is the Galois group of
$F[Y_i]/F[\zeta]$
and each $\sigma_k$ is as in the proof of Lemma 7, and $U$ is the Galois group of
$F[Y_i]/F[V_i]$. The subgroup $S$ of $G$ corresponding to
 $F[\sqrt[p]{N_i}]\subseteq F[V_i]$ in the Galois correspondence has index $p$ and
 contains~$U$. Therefore $S\supseteq H^p\rtimes U$, so $F[\sqrt[p]{N_i}]$ is contained in the fixed
 field of $H^p\rtimes U$, namely $F[W_i]$, where $W_i=\{\sqrt[p]{N_1},\dots,\sqrt[p]{N_\ell}\}\setminus \{\sqrt[p]{N_i}\}$.
 It follows from Lemma 1 that
 $N_1^{e_1}\cdots N_{\ell}^{e_\ell}\in F^p$, where all $0\leq e_j<p$ and $e_i=1$. By the rational root theorem, $F^p\cap D=D^p$,
 so $\S_p\cap D^p=\emptyset$ is violated.

If $m_i=2$ for at least one $i\in M^\sharp$, then (\ref{extension}) has been established. Likewise, if $p$ is odd, then (\ref{extension}) follows from (C). Suppose next that $4|m$ and $\S_2\cap (-D^2)=\emptyset$.
We claim that $\sqrt{-N_i}\notin F[V_i]$. If not, arguing as above, we see that $-N_1^{e_1}\cdots N_{\ell}^{e_\ell}\in F^2\subseteq K^2$, where all $0\leq e_j<2$ and $e_i=1$. On the other hand, $M\in K^2$ and $-1\in K^2$,
 so $-M\in K^2$ and therefore the product of all $N_j^{e_j+f_j}$, with $j\neq i$, must be in $K^2$. The uniqueness of $M$ in $\S_2\cap K^2$ forces $e_j=f_j$ for all
 $j\neq i$. Thus $-M\in F^2$ and hence $M\in\S_2\cap(-D^2)$, a contradiction. Thus (\ref{extension}) follows from (C) in this case as well.

Suppose finally that  $4|m_i$ for all $i\in M^\sharp$ and $M\in \S_2\cap (-D^2)$, say $M=-d^2$ with $d\in D$. Fix any $i\in M^\sharp$ and set $L_i=F[V_i]$.
It remains to decide when $N_i\in -4L_i^4$. Since $4|m_j$ for all $j\in M^\sharp$, the product of all $N_j^{f_j}$ with $j\neq i$
and $j\in M^\sharp$, belongs to $L_i^4$. Thus
$$N_i\in -4L_i^4\Leftrightarrow M\in -4L_i^4 \Leftrightarrow d^2\in 4L_i^4 \Leftrightarrow \pm 2d\in L_i^2,$$
and, by Lemma 1, this happens if and only if (\ref{spe}) holds.\qed

\section{Proofs of Theorems A and B}

\noindent{\sc Proof of Theorem A.} It is clear that (\ref{extension}) implies $\S_p\cap D^p=\emptyset$. Suppose
 $\S_p\cap D^p=\emptyset$ and let $K=F[\zeta]$, where $\zeta$ is a primitive $m$th root of unity.
By Lemmas 7 and 8, it suffices to show that $|\S_p\cap K^p|\leq 1$. Suppose not and let $M\neq N$ be in $\S_p\cap K^p$.
As $M,N$ have degree~$p$ over~$F$, we see that $p|[K:F]$. By Lemma 4, $\mathrm{Gal}(K/F)$ has a unique subgroup of index~$p$, so by the Galois correspondence, $K/F$ has a unique intermediate field of degree~$p$. We deduce
$F[\sqrt[p]{M}]=F[\sqrt[p]{N}]$ and Lemma 1 implies $MN^i\in F^p$ for some $i\in\Z$. Since $M\neq N$, this is disallowed by $\S_p\cap D^p=\emptyset$.\qed

\bigskip

\noindent{\sc Proof of Theorem B.} It is clear that $\S_2\cap D^2=\emptyset$ follows from (\ref{extension}).
Suppose $\S_2\cap D^2=\emptyset$. We will show that (\ref{extension}) holds if and only if $(N_1,\dots,N_\ell)$ is not 
2-defective.

By Lemmas 6, 7 and 8, we may restrict to the case when $|\S_2\cap K|=3$, in which case by Lemmas 5 and 6 there is a single
element $M\in \S_2\cap (-D^2)$, and we necessarily have $-1\notin F^2$ and $8|m$.

Now $-d^2=M=N_1^{f_1}\cdots N_\ell^{f_\ell}$, where $0\leq f_i<2$ and $M^\sharp=\{i\in J\,|\, f_i=1\}$ is nonempty.
Fix any $i\in M^\sharp$ and let $S^i_2$ stand for the analogue of $S_2$ corresponding to $\{N_1,\dots,N_\ell\}\setminus\{N_i\}$.
By the uniqueness of $M$ in $\S_2\cap (-D^2)$, we see that $\S_2^i\cap (-D^2)=\emptyset$. It follows from Lemma 6 that $|\S_2^i\cap K^2|\leq 1$.

Suppose first that $|\S_2^i\cap K^2|=0$. Set $V_i=\{\sqrt[m_1]{N_1},\dots,\sqrt[m_\ell]{N_\ell}\}\setminus \{\sqrt[m_i]{N_i}\}$,
and let $m[i]$ be the product of all $m_k$ such that $k\neq i$. Then $[K[V_i]:K]=m[i]$ by Lemma 7. Thus $F[V_i]$ is linearly disjoint
from $K$ over $F$. It follows that $\sqrt{-1}\notin F[V_i]$. For if $\sqrt{-1}\in F[V_i]$, then from $-1\notin F^2$ we deduce
that $1,\sqrt{-1}$ are $F$-linearly independent elements from $F[V_i]$, and hence $K$-linearly independent elements from $K[V_i]$,
which cannot be as $4|m$. Since $M\in -D^2$, we have $F[\sqrt{-1}]=F[\sqrt{M}]$. Thus $\sqrt{M}\notin F[V_i]$
and therefore $\sqrt{N_i}\notin F[V_i]$. If there is some $i\in M^\sharp$ such that $m_i=2$, this shows that (\ref{extension}) holds.
If, on the other hand, $4|m_i$ for all $i\in M^\sharp$, then (\ref{extension}) holds if and only if (\ref{spe}) fails, as in the proof of Lemma 8.

Suppose next that $|\S_2^i\cap K^2|=1$ and let $N\in S_2^i\cap K^2$. Note that $N\notin -D^2$. We have $N=N_1^{g_1}\cdots N_\ell^{g_\ell}$, where $g_i=0$, $0\leq g_j<2$ and $N^\sharp=\{j\in J\,|\, g_j=1\}$ is nonempty. Fix any $j\in N^\sharp$ and let
$S^{i,j}_2$ stand for the analogue of $S_2$ corresponding to $\{N_1,\dots,N_\ell\}\setminus\{N_i,N_j\}$. It is then clear that
$S_2^{i,j}\cap K^2=\emptyset$.  Set $V_{i,j}=\{\sqrt[m_1]{N_1},\dots,\sqrt[m_\ell]{N_\ell}\}\setminus \{\sqrt[m_i]{N_i},\sqrt[m_j]{N_j}\}$,
and let $m[i]$ (resp. $m[i,j]$) be the product of all $m_k$ such that $k\neq i$ (resp. $k\neq i,j$).
Then $[K[V_{i,j}]:K]=m[i,j]$ by Lemma~7. As above,
we deduce that $\sqrt{-1}\notin F[V_{i,j}]$. Since $4|m$ and $S^i_2\cap (-D^2)=\emptyset$, Lemma 8 ensures that $[F[V_i]:F]=m[i]$ as well as
$\sqrt{\pm N_j}\notin  F[V_{i,j}]$. We deduce from Lemma 2 that $\sqrt{-1}\notin F[V_i]$. The rest of the argument follows as in the above case.\qed

\bigskip


\section{Primitive elements}

Isaacs \cite{I} considered the problem of when $F[\a,\b]=F[\a+\b]$ for algebraic separable elements $\a,\b$ of degrees $m,n$ over $F$.
He proved that if $[F[\a,\b]:F]=mn$ (he actually assumed $\gcd(m,n)=1$ but used only the stated condition) but $F[\a,\b]\neq F[\a+\b]$
then $F$ has prime characteristic $p$ and the following conditions hold: $p|mn$ or $p<\min\{m,n\}$; if $m,n$ are prime powers, then $p|mn$; $p$ divides the order of the Galois group of a normal closure of $F[\a,\b]$.

The condition $p<\min\{m,n\}$ was later improved to $p<\min\{m,n\}/2$ by Divi\v{s} \cite{D}.

Using Isaacs' result we readily see that $F[\sqrt[m_1]{N_1},\dots,\sqrt[m_\ell]{N_\ell}]=F[b_1\sqrt[m_1]{N_1}+\cdots+b_\ell\sqrt[m_\ell]{N_\ell}]$ for any nonzero
$b_1,\dots,b_\ell\in F$ in Theorems A, B and C, provided the following conditions hold: $\chi(F)\neq p$ and $\S_p\cap D^p=\emptyset$ in
Theorem A; $\chi(F)\neq 2$, $\S_2\cap D^2=\emptyset$, and $(N_1,\dots,N_\ell)$ is not 2-defective in
Theorem B;   $\chi(F)\nmid m\varphi(m)$ (Euler's function), $\S_2\cap D^p=\emptyset$ for all $p\in\P$, and 
$(N_{i(2,1)},\dots,N_{i(2,\ell(2))})$ is not 2-defective in
Theorem C.

It {\em is}  possible that $F[\a,\b]/F$ be a finite Galois extension, $[F[\a,\b]:F]=[F[\a]:F][F[\b]:F]$ and still $F[\a,\b]\neq F[\a+\b]$.
A family of examples can be found in \cite[Example 2.3]{CS}.


\end{document}